\documentclass[letterpaper,12pt,leqno]{article}

\usepackage{fullpage}

\usepackage{hyperref}

\usepackage{amsmath}
\usepackage{amsfonts}
\usepackage{amssymb}
\usepackage{tabularx}
\usepackage{color}
\usepackage{fancyhdr}
\usepackage[all]{xy}
\usepackage{enumerate}

\def\theoremparent{section}
\usepackage{notes}

\numberwithin{equation}{subsection}
\numberwithin{theorem}{subsection}

\def\A{\mathbb{A}}
\def\abs{\mathrm{abs}}
\def\AbsAc{\mathrm{AbsAc}}
\def\Acyc{\mathrm{Acyc}}
\def\b{\mathrm{b}}
\def\bb{\mathfrak{b}}

\def\Bar{\mathrm{Bar}}
\def\Ber{\mathrm{Ber}}

\def\cD{\mathcal{D}}
\def\cg{\mathrm{cg}}
\def\CMod{\textrm{-}\mathrm{CMod}}
\def\cnt{\mathrm{cnt}}
\def\CntAc{\mathrm{CntAc}}
\def\co{\mathrm{co}}
\def\CoAc{\mathrm{CoAc}}
\def\Cob{\mathrm{Cob}}

\def\coRees{\mathrm{coRees}}

\def\d{\partial}
\def\D{\mathrm{D}}
\def\Der{\mathrm{Der}}
\def\dgVect{\mathrm{dgVect}}
\def\dgsVect{\mathrm{dgsVect}}

\def\Dmon{D_{\mathrm{mon}}}
\def\dR{\mathrm{dR}}
\def\ev{\mathrm{ev}}
\def\fd{\mathrm{fd}}
\def\fg{\mathrm{fg}}
\def\g{\mathfrak{g}}

\def\GG{\mathbb{G}}
\def\gr{\mathrm{gr}}
\def\grco{{\gr\textrm{-}\co}}
\def\Htp{\mathrm{Htp}}
\def\Ind{\mathrm{Ind}\textrm{-}}
\def\Inj{\mathrm{Inj}}

\def\MC{\textrm{-}\mathrm{MC}}
\def\nn{\mathfrak{n}}
\def\odd{\mathrm{odd}}

\def\op{\mathrm{op}}

\def\Perf{\mathrm{Perf}}

\def\proof{\noindent{\em Proof:}\ }

\def\qed{\hfill\lower 1em\hbox{$\square$}\vskip 1em}
\def\Rees{\mathrm{Rees}}

\def\s{\mathfrak{s}}

\def\shp{\#\!}
\def\shpInj{\textrm{-}\#\Inj}
\def\Sing{\mathrm{Sing}}

\def\str{\mathrm{str}}
\def\sVect{\mathrm{sVect}}

\def\tr{\mathrm{tr}}
\def\toto{\text{\ \raise0.2em\hbox to 0pt {$\to$}\lower0.2em\hbox{$\to$}\ }}

\def\ss#1{\!\<#1\>}

\begin{document}

\title{Monadicity of localization for Lie super-algebras $\gl(m, n)$}
\author{Slava Pimenov}
\date{\today}
\titlepage
\maketitle

\tableofcontents

\vskip 5em
\setcounter{section}{-1}
\section{Introduction}
One version of Beilinson-Bernstein localization theory for a semisimple algebraic group $G$ can be stated as follows
(\cite{BK}, \cite{BMR}, \cite{BN}).
Let $B$ be a Borel subgroup and $N$ its unipotent radical. Denote by $\tilde X = G/N$ the affine flag variety
and $\Dmon(\tilde X)$ the category of monodromic $D$-modules on $\tilde X$. Functors of global sections and localization
between categories of $U\g$-modules and monodromic $D$-modules form a Frobenius adjunction, in the sense that they are
adjoint to each other from both sides up to twist by the dualizing sheaf $\omega_{\tilde X}$.

The composition of localization and global sections functors give both a monad and a comonad $W$ on $\Dmon(\tilde X)$, and
the Barr-Beck theorem identifies the derived category $D(U\g\Mod)$ with the categories of $W$-modules and $W$-comodules
in $D(\Dmon(\tilde X))$. In particular for a regular weight $\lambda$ the specialization of the monad $W_\lambda$ can be
seen to be the group algebra of the Weyl group of $G$, its action identifies monodromic $D$-modules with the monodromies on
the $W$-orbit of $\lambda$. This recovers the usual equivalence of the twisted $D$-modules with regular dominant weight
and representations of $\g$ with corresponding central character.

Moreover, identifying the category of $G$-equivariant endofunctors on $\Dmon(\tilde X)$ with the $D$-modules on the double
quotient $N\backslash G/N$ the Weyl monad $W$ corresponds to the differential operators on $N\backslash G/N$. The orbit
structure of the double quotient is used to obtain information regarding the monad $W$ itself.

In case of the super-groups we no longer have the Frobenius adjunction, however we still can define one of the adjoints
(left of right depending on the setting) after replacing the derived categories of representations and monodromic $D$-modules
with suitable Ind-completions. The main result of this paper is the ind-monadicity of the right localization functor, which
allows us to identify the ind-completion of the bounded derived category of $U\g$-modules with the ind-completion of
$W$-modules in $\Dmon(\tilde X)$ (\ref{thm_ind_monadic}). We consider here only Lie super-algebras of type $A$, however we expect the theorem to
hold in other types as well, since the main point of the proof is the reduction to the classical group underlying the super-group.

In order to identify the Weyl monad with the differential operators on the double quotient, we believe we need to work with
filtered $D$-modules or even better Hodge modules. We will investigate this question elsewhere.

\vskip 1em
\begin{nparagraph}[Outline of the paper.]
We start by recalling the definition and basic properties of coderived categories (\cite{Pos}), as well as a Koszul duality
theorem where these categories arise. Next investigate relation between $D$-modules and $\Omega$-modules on super-manifolds
and establish equivalences between appropriate derived categories analogous to those shown by Positselski. In particular
we show that the derived category of $D$-modules is equivalent to the coderived category of $\Omega$-modules (\ref{thm_D_Omega})
and the equivalence of coderived categories of filtered $D$-modules and filtered $\Omega$-modules (\ref{thm_D_Omega_filt}).

We finish the section with discussion of equivariant $D$-modules, and illustrate why we believe that filtration is necessary
to obtain the sound theory of $D$-modules on the quotient stacks.

In the second section we show an equivalence between coderived categories of $U\g$-modules and modules over the cochain algebra
$C^\bullet(\g)$ for a Lie super-algebra $\g$. Then we define localization functors both on the $D$-module and $\Omega$-module sides
and show that they are compatible with these equivalences. Finally we define the notion of ind-monadicity of a functor and show
that the right localization functor in monodromic $D$-modules is ind-monadic.

In conclusion we illustrate the main theorem in the case of $\sl(1, 1)$ where we give an explicit description of the Weyl monad
as an algebra in monodromic $D$-bimodules.

\end{nparagraph}

The author would like to thank Kobi Kremnitzer for the continued interest and support of this work.

\vfill\eject
\section{$D$-modules and $\Omega$-modules}

\subsection{Notations and setting}
Throughout this paper we work over an algebraically closed field $k$ of characteristic $0$.

\begin{nparagraph}[Category $\sVect$.]
By a super vector space we mean a $\Z$-graded vector vector space $V = \bigoplus_i V_{s=i}$, where $V_{s=i}$ or simply $V_i$,
when there is no ambiguity what grading is being used, is the component of super-degree $i$. We will write $V\ss{n}$ for the
super vector space with the $n$-shifted grading, i.e.
$$
V\ss{n}_i = V_{i+n}.
$$
By $\sVect$ we denote the category of super vector spaces and $k$-linear maps of degree $0$. The degree of an element $v \in V$
will be denoted by $|v|_s$ or simply $|v|$, when the implied grading is clear.

We have the standard symmetric monoidal structure on $\sVect$:
$$
(V \tensor W)_n = \bigoplus_{i+j=n} V_i \tensor W_j,
$$
with the symmetry map $V \tensor W \to W \tensor V$ given by the Koszul sign rule. Namely, for $v \in V$ and $w \in W$
$$
v \tensor w = (-1)^{|v|_s |w|_s} w \tensor v.
$$

\end{nparagraph}

\begin{nparagraph}[Category $\dgsVect$.]
\label{cat_dgsVect}
Let $V^\bullet$ be a complex of super vector spaces, with cohomological grading, i.e. with differentials $d\from V^i \to V^{i+1}$.
The cohomological degree of an element $v \in V^\bullet$ will be denoted by $\bar v$. As customary, the cohomological shift
will be denoted by
$$
V[n]^i = V^{n+i}.
$$

We will write $\dgsVect$ for the dg-category of complexes of super vector spaces. We have the standard symmetric monoidal structure
$$
(V \tensor W)^n = \bigoplus_{i+j=n} V^i \tensor W^j,
$$
and the symmetry map is given by the Koszul sign with respect to the total grading:
$$
v \tensor w = (-1)^{(\bar v + |v|_s) (\bar w + |w|_s)} w \tensor v.
$$

\end{nparagraph}

\begin{nparagraph}[Filtered complexes and Rees construction.]
\label{rees_functor}
Consider a complex $V^\bullet$ equipped with an increasing cocomplete filtration by subcomplexes $F_p V^\bullet$. Rees construction
is the canonical functor from the quasi-abelian category of filtered complexes to its abelian envelop. Explicitly, it is given
by the graded object
$$
\Rees (V^\bullet, F) = \bigoplus_p F_p V^\bullet.
$$

The degree of $v \in \Rees(V^\bullet, F)$ coming from the Rees construction will be denoted by $|v|_h$. We can think of 
$\Rees (V^\bullet, F)$ as a $k[h]$-module with multiplication by $h$ given by inclusion $F_p V \into F_{p+1} V$.
Here $h$ is a formal variable of degrees
$$
\bar h = |h|_s = 0,\quad |h|_h = 1.
$$
The subscript $h$ in the notation for the degree, refers both to the formal variable $h$ as well as to the Hodge filtration that
we will be concerned with. The essential image of the Rees functor is formed by flat $k[h]$-modules, and because of the
cocompleteness assumption it is fully faithful. Indeed, given a flat $k[h]$-module $M$ we can recover $V^\bullet$ as $V^\bullet = \colim M_i$.

Let $V$ be a filtered complex, then we will denote the complex with shifted filtration (or Tate twist) by $V(n)$:
$$
F_p V(n) = F_{p + n} V.
$$

For a filtered algebra $A$ in $\dgsVect$ the Rees construction $\Rees (A, F)$ is a flat central $k[h]$-algebra. If $M$ is a
filtered $A$-module with a compatible filtration, i.e. $F_p A \tensor F_q M \to F_{p+q} M$, then $\Rees (M, F)$ is a
$\Rees (A, F)$-module.

\end{nparagraph}

\vskip 3em
\subsection{Enhancements of the Derived Category}

Let $A$ be an algebra in $\dgsVect$. We will denote by $\Htp(A\Mod)$ the homotopy category, i.e. the localization of $A$-modules
with respect to homotopy equivalences, and by $\D(A\Mod)$ the unbounded derived category, i.e. the localization of $\Htp(A\Mod)$
with respect to quasi-isomorphisms, or equivalently the quotient of $\Htp(A\Mod)$ by the thick subcategory of acyclic objects
$\Acyc(A\Mod)$. We will now recall the definitions of coderived, contraderived and absolutely derived
categories.

\begin{nparagraph}[Absolutely acyclic, coacyclic and contra-acyclic modules.]
The category of absolutely acyclic $A$-modules is the minimal subcategory of $A\Mod$, containing total complexes of short exact
sequences of $A$-modules and closed under shifts, cones and retracts. In particular it contains total complexes of all finite
exact sequences. We will denote this category as $\AbsAc(A\Mod)$.

The category of coacyclic modules is the minimal subcategory of $A\Mod$ containing absolutely acyclic modules and closed under
all operations listed above and under taking filtered colimits. It will be denoted by $\CoAc(A\Mod)$.

Similarly, the contra-acyclic modules is the minimal subcategory containing absolutely acyclic modules and closed under shifts,
cones, retracts and filtered limits. It will be denoted by $\CntAc(A\Mod)$.

We will say that a morphism of $A$-modules and absolutely acyclic, coacyclic or contra-acyclic if its cone belongs to the
corresponding subcategory.

\end{nparagraph}

\begin{nparagraph}[Absolute derived, coderived and contraderived categories.]
The absolute derived, coderived and contraderived categories are defined as quotients of the homotopy category
by the respective subcategories of acyclics.
\begin{align*}
\D^\abs(A\Mod) &= \Htp(A\Mod) / \AbsAc(A\Mod)\\
\D^\co(A\Mod) &= \Htp(A\Mod) / \CoAc(A\Mod)\\
\D^\cnt(A\Mod) &= \Htp(A\Mod) / \CntAc(A\Mod)
\end{align*}

From the definitions it is clear we have the following diagram of localizations:
$$
\xymatrix{
&\D^\abs(A\Mod) \ar[dl] \ar[dr] &\\
\D^\co(A\Mod) \ar[dr] &&\D^\cnt(A\Mod) \ar[dl]\\
&\D(A\Mod)&
}
$$
In general none of these functors are equivalences, however there are many results relating these categories in special cases.
We will recall here some of the results that will be relevant for us.

For an algebra $A$ in $\dgsVect$ we will denote $\shp A$ the graded algebra obtain from $A$ by forgetting the differential.
The full subcategory of $A\Mod$ formed by modules $M$, such that $\shp M$ is an injective $\shp A$-module will be
denoted by $A\shpInj$.

We say that $\shp A$ is {\em Gorenstein} if the classes of $\shp A$-modules of finite projective and finite injective dimensions
coincide. Denote by $A\Mod_\fd$ the full subcategory of $A\Mod$ formed by modules $M$ such that $\shp M$ is of finite projective and injective dimensions
as an $\shp A$-module.

\end{nparagraph}

\begin{proposition}
\label{prop_posit}
\begin{enumerate}[a)]
\item Let $A$ be such that $\shp A$ is of finite cohomological dimension, then subcategories $\AbsAc(A\Mod)$, $\CoAc(A\Mod)$ and
$\CntAc(A\Mod)$ coincide.
\vskip 1em
\item Let $A$ be such that $\shp A$ is Noetherian, then we have the semi-orthogonal decomposition
$$
\Htp(A\Mod) = \left< \Htp(A\shpInj), \CoAc(A\Mod)\right>,
$$
in the sense that $\Hom(K, I)$ is an acyclic object in $\dgVect$ for any $I \in A\shpInj$ and $K \in \CoAc(A\Mod)$.

Moreover, the following composition of functors is an equivalence:
$$
\Htp(A\shpInj) \to \Htp(A\Mod) \to \D^\co(A\Mod).
$$

\vskip 1em
\item\label{test} Let $A$ be such that $\shp A$ is Gorenstein, then the following functors are equivalences:
$$
\xymatrix{
\D^\co(A\Mod) & \D^\abs(A\Mod_\fd) \ar_\equ[l] \ar^\equ[r] & \D^\cnt(A\Mod).
}
$$
\end{enumerate}
\end{proposition}

We would like to point out that equivalences of coderived and contraderived categories in (\ref{test}) do not imply that
$\CoAc(A\Mod)$ and $\CntAc(A\Mod)$ coincide as subcategories of $\Htp(A\Mod)$.

Let us also recall the compact generation of the coderived category.

\begin{proposition}
\label{prop_compact_gen}
Let $A$ be such that $\shp A$ is Noetherian, then the coderived category $\D^\co(A\Mod)$ is compactly generated by the
absolute derived category of finitely generated $A$-modules $\D^\abs(A\Mod_\fg)$.
\end{proposition}

In particular, when $A$ is concentrated in cohomological degree $0$, we have $\D^\co(A\Mod) \isom \Ind\D^\b(A\Mod_\fg)$,
since an acyclic complex in $\D^\b(A\Mod)$ is a finite exact sequence, and hence absolutely acyclic.

\begin{nparagraph}[Singularity category.]
Recall that the unbounded derived category $\D(A\Mod)$ is compactly generated by the perfect complexes $\Perf(A\Mod)$.
When $A$ is concentrated in cohomological degree $0$, the singularity category of $A$ is defined as the quotient
$$
\Sing(A\Mod) = \D^\b(A\Mod_\fg) / \Perf(A\Mod).
$$
In other words, in this case the difference between derived and coderived categories is the $\Ind$completion of the singularity
category.

\end{nparagraph}

\begin{nparagraph}[Koszul duality between modules and comodules.]
There are many different versions of Koszul duality, here we will recall the duality between the categories of comodules over
a conilpotent coaugmented coalgebra $(C, \Delta)$ in $\dgsVect$ and modules over the reduced cobar-construction $A = \wbar{\Cob(C)}$.
Let $\epsilon\from k \to C$
denote the coaugmentation map, then we say that $C$ is {\em conilpotent} if the composition of iterated comultiplication map and
projection to the quotient
$$
\xymatrix{
C \ar^-{\Delta}[r] & C^{\tensor N} \ar@{>>}[r] & (C/k)^{\tensor N}
}
$$
vanishes for large $N$. The reduced bar-complex is defined as
$$
A\ =\ \wbar{\Cob(C)}\ =\ \bigoplus_{n}\ ((C/k)[-1])^{\tensor n},
$$
with the differential induced by comultiplication in $C$, and multiplication $\mu$ given by concatenation.

Consider complex $C \tensor A$ with the differential defined as composition
$$
\xymatrix@C=3em{
d\from C \tensor A \ar^-{\Delta\tensor\id}[r] & C \tensor C \tensor A \ar[r] & C \tensor (C/k)[1] \tensor A \ar^-{\id\tensor\mu}[r] & C \tensor A[1].
}
$$
This complex is a right $A$-modules and left $C$-comodule. Similarly, we define complex $A \tensor C$ in left $A$-modules and right $C$-comodules.

\begin{proposition}
\label{prop_KD_bar}
Let $C$ and $A$ be as above, then the pair of adjoint functors
\begin{align*}
(C \tensor A) \tensor_{A} -\ &\from \Htp(A\Mod) \to \Htp(C\Comod)\\
(A \tensor C) \tensor_{C}^{\co} -\ &\from \Htp(C\Comod) \to \Htp(A\Mod)
\end{align*}
induce equivalence of categories $\D^\co(C\Comod)$ and $\D(A\Mod)$. Here $\tensor_C^\co$ denotes the cotensor product of $C$-comodules.

\end{proposition}

\end{nparagraph}

\begin{nparagraph}[Duality between filtrations and bicomplexes.]
\label{filt_bicomplex}
As a simple but useful example to illustrate the difference between derived and coderived categories and provide some intuition to the reader
unfamiliar with these notions we apply the Koszul duality above to the case of filtered complexes and bicomplexes.

Consider coalgebra $C = k[\epsilon]$, with $\bar\epsilon = -1$ and $|\epsilon|_h = 1$, and comultiplication
$\Delta(\epsilon) = 1 \tensor \epsilon + \epsilon \tensor 1$. Denote the $k$-dual algebra by $B = C^*$. Since $C$ is of finite dimension,
categories of $C$-comodules and $B$-modules are equivalent.

We also consider the algebra $R = k[h]$, with $\bar h = 0$ and $|h|_h = 1$, and the graded dual coalgebra $S = \Hom_\gr(R, k)$.
It is clear that $R$ is isomorphic to the reduced bar-complex $\wbar{\Bar(C)}$ and $B$ is quasi-isomorphic to $\wbar{\Bar(S)}$.

The category of $C$-comodules (or equivalently $B$-modules) can be identified with the category of bicomplexes, say with horizontal
differential $d$ and vertical differential induced by the coaction of $\epsilon$. The cohomological grading then corresponds to the
total grading of the bicomplex.

As we mentioned in paragraph \ref{rees_functor}, $R$-modules contain a full subcategory equivalent to the category of complexes
with cocomplete filtration. Now, let us look at $S$-comodules. Again we start with a filtered complex $(V, F)$ and form the
graded complex
$$
\coRees (V, F)\ =\ \bigoplus_i\ V/F_i V,
$$
with the term $V/F_i V$ placed in $h$-degree $i$. The coaction of $S$ is given by
$$
\xymatrix{
v \ar@{|->}^-a [r] & \sum_j\ (h^j)^* \tensor [v]_{i+j},
}
$$
where $v \in V/F_i V$ and $[v]_{i+j}$ is the class of $v$ in $V/F_{i+j}V$. The cocompleteness of filtration ensures that this sum
is finite.

Given a coflat $S$-comodule $N$, we can recover complex with a complete and cocomplete filtration by taking inverse limit $V = \lim_i N_i$.

\vskip 1em

Let us now describe the adjunction of proposition \ref{prop_KD_bar} in this case. For simplicity consider flat $R$-module corresponding
to a filtered complex $(V, F)$ concentrated in a single cohomological degree. Applying functor $(C \tensor R) \tensor_R -$, we obtain the following bicomplex.
$$
\xymatrix{
... & F_i && \\
& \epsilon F_i \ar^-d[r] \ar@{-->}[u] & F_{i+1} & \\
&& \epsilon F_{i+1} \ar^-d[r] \ar@{-->}[u]& F_{i+2} \\
&&& ...
}
$$
Applying the adjoint functor we obtain a flat $R$-module given by the totalization of the bicomplex above with the filtration by
horizontal slices. We have a similar picture for $S$-comodules and $B$-modules. Notice that since $S$-comodules have finite cohomological
dimension we have $\D^\co(S\Comod) \equ \D(S\Comod)$. Putting it together:
\vskip 1em
$$
\xymatrix@C=4em@M=0.5em@R=5em{
\displaystyle
\left\{ \vcenter{\hbox to 12em {\hfil complexes with\hfil}\hbox to 12em {\hfil complete and cocomplete\hfil}\hbox to 12em{\hfil filtrations\hfil }}\right\}
\ar@{^{(}->}[r] \ar@{_{(}->}[d] & \D(S\Comod) \ar^\equ@<0.2em>[r] \ar@{_{(}->}[d] & \ar@<0.2em>[l] \D(B\Mod) \ar@{_{(}->}[d] \\
\left\{ \vcenter{\hbox to 8em {\hfil complexes with\hfil}\hbox to 8em {\hfil cocomplete\hfil}\hbox to 8em{\hfil filtrations\hfil }}\right\}
\ar@{^{(}->}[r] & \D(R\Mod) \ar^\equ@<0.2em>[r] & \ar@<0.2em>[l] \D^\co(B\Mod).
}
$$

\end{nparagraph}
\vskip 2em

\begin{nparagraph}[Relative cobar construction.]
Let $C$ be a commutative Hopf algebra, then the tensor product $\tensor_k$ can be thought of as a symmetric monoidal structure in right $C$-comodules,
with the coaction on the product given as the composition
$$
\xymatrix{
V \tensor W \ar[r] & (V \tensor C) \tensor (W \tensor C) \ar^-{\mu_C}[r] & (V \tensor W) \tensor C.
}
$$

Let $A$ be an algebra in right $C$-comodules with respect to this monoidal structure with the coaction $\alpha\from A \to A \tensor C$.
Explicitly this means that the following diagram commutes.
$$
\xymatrix@C=4em{
A \tensor A \ar^-{\mu_A}[r] \ar_{\alpha \tensor \alpha}[d] & A \ar^{\alpha}[d] \\
(A \tensor C) \tensor (A \tensor C) \ar^-{\mu_A \tensor \mu_C}[r] & A \tensor C.
}
$$
The category of (either left or right) $A$-modules in right $C$-comodules with respect to this monoidal structure will be denoted $(A, C)\MC$.

The Koszul duality between right $C$-comodules and right $B = \wbar{\Cob(C)}$-modules sends $A$ to a $B$-module
$$
\tilde A\ =\ \bigoplus_{n}\ A \tensor (C/k)^{\tensor n},
$$
with the differential induced by the right coaction of $C$ on $A$ and the comultiplication in $C$. In fact this is a dg-algebra over $B$.
To see this it will be convenient to work with an alternative construction of $\tilde A$.

Consider an algebra $H = A \tensor C$. It is an $A$-bimodule with the obvious left multiplication and the right multiplication given by
$$
c \tensor a \mapsto \alpha(a)_0 \tensor \alpha(a)_1 c,
$$
here we use Sweedler notation for the coaction. The counit map $C \to k$ induces the map of $A$-bimodules $H \to A$, and we denote
by $H^+$ its kernel. $H$ is a coalgebra in $A$-bimodules with the comultiplication $H = A \tensor C \to H \tensor_A H = A \tensor C \tensor C$
defined as
$$
\Delta_H(a \tensor c) = a \tensor \Delta_C(c).
$$

The differential $d$ in the cobar complex $\Cob_A(H) = \bigoplus_n H \tensor_A \ldots \tensor_A H$ is induced by the comultiplication $\Delta_H$.
Since the unit map $k \to C$ is a map of coalgebras the element $1_C \in C[-1] \subset \Cob_A(C)$ is a Maurer-Cartan element and we can
deform the differential in the cobar complex by setting $d' = d + [1_C, -]$.

Identifying $H^+$ with the quotient $H / A$ of left $A$-modules we find that both the cobar differential $d$ and the deformed differential
$d'$ induce a well defined differential on the reduced cobar complex
$$
\wbar{\Cob_A(H)} = \bigoplus_n H^+ \tensor_A \ldots \tensor_A H^+,
$$
making it into a dg-algebra. It is immediate to see that this complex with the differential induced by $d'$ is isomorphic to $\tilde A$ defined before.

\end{nparagraph}

\begin{nparagraph}[Duality for modules over filtered algebras.]
\label{duality_filt_A}
Now we would like to upgrade the discussion of section \ref{filt_bicomplex} to the case of modules over algebras. First notice that $C = k[\epsilon]$
is a commutative Hopf algebra, so we are in the setting of the previous section.
Moreover, since $C$ is cocommutative we will not distinguish between left and right $C$-comodules.

Let $A$ be an algebra in $C$-comodules, then
using construction above section we obtain an $R = \wbar{\Cob(C)}$-algebra $\tilde A = R \tensor A$. Observe that in this particular
case this $R$-algebra is central. Furthermore, for a module $M \in (A, C)\MC$ we obtain an $\tilde A$-module $\tilde M = R \tensor M$.
We want to compare categories $(A, C)\MC$ and $\tilde A\Mod$.

Using the natural map $A \to C \tensor \tilde A = C \tensor R \tensor A$, we see that for an $\tilde A$-module $M$ the tensor product
$C \tensor M$ is an object of $(A, C)\MC$.

\end{nparagraph}

\begin{proposition}
\label{prop_filt_A}
The adjoint pair of functors
\begin{align*}
(C \tensor R) \tensor_{R} -\ &\from \Htp(\tilde A\Mod) \to \Htp((A, C)\MC)\\
(A \tensor C) \tensor_{C}^{\co} -\ &\from \Htp((A, C)\MC) \to \Htp(\tilde A\Mod)
\end{align*}
induce equivalence of categories $\D^\co(\tilde A\Mod)$ and $\D^\co((A, C)\MC)$.
\end{proposition}
\proof
The proof of proposition \ref{prop_KD_bar} applies here without change, since all involved constructions respect the action of $A$.
We will provide it here for convenience of the reader.

Since filtered colimits in the categories of $A$-modules and $\tilde A$-modules are calculated at the level of underlying vector spaces
we see that both functors preserve filtered colimits and hence send coacyclic objects to coacyclic.

We need to show that the unit and counit maps for the adjunction are coacyclic. Consider the map $M \to C \tensor R \tensor M$ in
$(A, C)\MC$. The cone of this map is the direct sum totalization of the bicomplex
$$
B^{pq}\ =\ \begin{cases}
(C \tensor h^p) \tensor M^q,&\text{if $p \ge 0$}\\
M^q,&\text{if $p = -1$},
\end{cases}
$$
with vertical differential induced from $M$ and horizontal from $C \tensor R$, namely $d(\epsilon h^p) = h^{p+1}$. Consider the canonical
filtration of $B$ by vertical slices. Each term of the filtration is a finite exact sequence of objects in $(A, C)\MC$, hence
its totalization is an absolutely acyclic
object of $(A, C)\MC$, and therefore the cone is absolutely acyclic.

Now take $N \in \tilde A\Mod$ and consider the map $R \tensor C \tensor N \to N$, we will show that it is absolutely acyclic.
Look at the following two-stage filtration $G_\bullet K$ of the cone $K$ of this map
$$
G_0 = R \tensor \epsilon \tensor N \into R \tensor C \tensor N \into K = G_2.
$$
The terms of this filtration are closed under the action of $\tilde A$, but not under differential. However, $d(G_i) \subset G_{i+1}$,
so that it induces a differential on the associated graded
\begin{equation}
\label{shexseq_1}
R \tensor N \to R \tensor N \to N,
\end{equation}
which is a resolution of $N$ as an $R$-module. Let us modify filtration $G$ in the following way, set $\tilde G_0  = G_0 + d(G_0)$ and
$\tilde G_1 = K$. Since the complex in (\ref{shexseq_1}) is exact we see that both $\tilde G_0$ and $K / \tilde G_0$ are contractible,
therefore the cone $K$ is absolutely acyclic.

\qed

\begin{remark}[$S^1$-equivariant algebraic geometry.]
We can think of the Hopf algebra $C$ as the cohomology ring of topological circle $H^\bullet(S^1)$. In this way one can think
of vector spaces and algebras with coaction of $C$ as $S^1$-equivariant objects (see for instance \cite{TV}). In light of the
previous proposition, working in the $S^1$-equivariant setting in equivalent to working over the filtration stack $\A^1 / \GG_m$.

\end{remark}

\vskip 3em
\subsection{$D$- and $\Omega$-modules}
\label{ss_D_Omega}

We start discussion of $D$-modules by recalling the basic facts in the non-derived case (\cite{Penkov1}). A smooth $\Z$-super-manifold is a pair
$(X, E)$ consisting of a smooth classical manifold $X$, equipped with an even $\Z$-grading on the structure sheaf $\O_X$, and
a bundle $E$ of finite rank over $X$ concentrated in odd degrees. We will call the sheaf of algebras
$A = \Sym^\bullet_{\O_X} E$ (with the total degree sign convention as specified in \ref{cat_dgsVect}) the sheaf of functions
on the super-manifold $(X, E)$.

Let $X$ be affine, denote $\Der(A) \subset \End(A)$ the space of derivations, it is easy to see that
$$
\Der(A) = A \tensor_{\O_X} \left(\Der(\O_X) \oplus E\right).
$$
The algebra of differential operators is defined as in the classical case. For instance, consider subspace $A \oplus \Der(A) \subset \End(A)$,
where the first summand corresponds to left multiplication. It is an $A$-subbimodule and a Lie subalgebra, and so we can define
$$
D_A = \left(\bigoplus_n (A \oplus \Der(A))^{\tensor_A n}\right)\ \bigg/ \ \left(uv - vu = [u, v],\ 1_{A} = 1\right).
$$

Localizing in $X$ if necessary, we may assume that $E$ is a trivial bundle with fiber $V \in \sVect$, in which case we have
the decomposition of differential
operators $D_A = D_X \tensor D_V$. The second factor $D_V$ is the algebra of differential operators on the odd super vector space $V$.
It is straightforward to see that $D_V = \End_k(\Sym^\bullet V)$, the space of $k$-algebra endomorphisms, therefore by Morita equivalence we find that $D_V\Mod \equ \sVect$,
hence the category $D_A\Mod$ is equivalent to the category of $D_X$-modules in $\sVect$.

\begin{nparagraph}[De Rham and Spencer complexes.]
Consider the following two graded super-vector spaces
\begin{align*}
Sp\ &=\ \Sym^\bullet_A\ (\Der(A)[1]),\\
\Omega\ &=\ \Sym^\bullet_A\ (\Omega^1[-1]).
\end{align*}

For any left $D$-module $M$ we can form the De Rham complex $\Omega(M) = \Omega \tensor_A M$ with the differential induced by
the de Rham differential in $\Omega$, differential in $M$ and the natural map $A \to \Omega^1 \tensor \Der(A)$. It is a left
$\Omega$-module.
For a right $D$-module $N$ we can form the Spencer complex $Sp(N) = \Hom_{D^\op}(\Omega(D), N)$, which is a right $\Omega$-module.

The special important case of these complexes is when they are applied to the $D$-bimodule $D$. In this case
the Spencer complex $Sp(D)$ is quasi-isomorphic to the algebra $A$ as a right $\Omega$-module and the De Rham complex
$\Omega(D)$ is quasi-isomorphic to the dualizing sheaf $\omega$, i.e. the shifted Berezinian line bundle $\Ber_A[-\dim X]$.

Let us point out that in general $Sp(N)$ is the completion of the complex $N \tensor_A Sp$ with respect to the stupid
filtration on $N$. In particular when $N$ is bounded from above or $A$ is purely even we have $Sp(N) \isom N \tensor_A Sp$.

Algebra $D_A$ and the space $Sp$ are equipped with increasing filtrations induced by tensor powers, and the space $\Omega$
with a decreasing filtration. They induce filtrations by $\Omega$-submodules on $Sp(D)$ and $\Omega(D)$:
\begin{align*}
F_n Sp(D_A)\ &=\ \sum_{j+i = n}\ \Sym^i_A (\Der(A)[1]) \tensor F_j D_A,\\
F_n \Omega(D_A)\ &=\ \sum_{j - i = n}\ \Sym^i (\Omega^1[-1]) \tensor F_j D_A.
\end{align*}

We will need the following lemma.

\end{nparagraph}

\begin{lemma}
\label{spencer_coacyclic}
The natural map $A \to Sp(D)$ is coacyclic in $\Omega$-modules.
\end{lemma}
\proof
The proof is the same as in the classical case. The above filtration on the Spencer complex is a bounded from below filtration
by $\Omega$ submodules. So the kernel $K$ of the map in question $Sp(D) \to A$ is a filtered colimit of $\Omega$-submodules $F_n K$. It is enough
to show that each $F_n K$ is coacyclic and for that it is enough to show that each quotient $F_n / F_{n-1}$ is coacyclic.
The latter is a finite exact sequence of $A$-modules, hence it is absolutely acyclic.

\qed

\vskip 1em
Positselski showed that in the classical case we have equivalence of the derived category of $D$-modules and coderived category of
$\Omega$-modules on a smooth scheme. The same holds in the super case.

\begin{theorem}
\label{thm_D_Omega}
Let $A$ be the ring of functions on a smooth $\Z$-super-manifold, then we have an equivalence
$$
\D(D_A\Mod) \equ \D^\co(\Omega_A\Mod).
$$
\end{theorem}
\proof
Localizing in $X$ if necessary we may assume that the odd bundle is trivial with the fiber $V$, so that $A \isom \O_X \tensor \Sym^\bullet V$.
We will separate even and odd variables and establish the equivalence in question as the composition of the following two equivalences
$$
\xymatrix@C=6em{
\D(D_A\Mod) \ar^-{\Hom(\Omega_X, -)}@<0.2em>[r] & \ar^-{D_X \tensor -}@<0.2em>[l] \D^\co(\Omega_X \tensor D_V\Mod)
\ar^-{\Omega_V \tensor -}@<0.2em>[r] & \ar^-{\Hom(D_V, -)}@<0.2em>[l] \D^\co(\Omega_A\Mod).
}
$$

\vskip 1em
\textit{Even case.}
This is the equivalence established by Positselski, we provide the proof here to illustrate the difference with the odd case.

Since cohomological dimension of $D_A$ is finite and $D_A$ is concentrated in cohomological degree $0$ we have $\Acyc(D_A\Mod) \equ \CoAc(D_A\Mod)$.
As was mentioned before, $\Hom(\Omega_X, N) \isom Sp \tensor N$, therefore both adjoint functors commute with filtered colimits
and hence send coacyclic objects to coacyclic.

In order to show that the natural map $D \tensor Sp \tensor N \to N$ is acyclic, use the increasing filtration on $D\tensor Sp$. The quotients
of this filtration are acyclic except the $0$-term which is isomorphic to $A$. Using the spectral sequence argument we find that the
counit map is a quasi-isomorphism.

Now we need to show that the unit $M \to Sp \tensor D \tensor M$ is coacyclic in $\Omega_X \tensor D_V$-modules. First, filter
$M$ by submodules $G_p M = (\Omega^{\ge p}_X \tensor D_V) \cdot M$. This is a finite filtration, so it is enough to establish coacyclicity
for the associated graded module $\gr^G M$. Each quotient of the filtration $\gr^G_i M$ is killed by $\Omega^{\ge 1}_X$, and therefore it is
an $\O_X \tensor D_V$-module. The unit map $\gr^G_i M \to Sp \tensor D \tensor \gr^G_i M$ is the tensor product over $\O_X$ of $\gr^G_i M$ and the
coacyclic map $\O_X \to Sp(D_X)$ (lemma \ref{spencer_coacyclic}), hence is also coacyclic.

\vskip 1em
\textit{Odd case.}
Now we move to the second adjunction. Since the algebra $D_V$ is finite-dimensional both functors of the adjunction again preserve
filtered colimits and therefore preserve coacyclics.

Denote by $D^*_V = \Hom_{\O_V}(D_V, \O_V)$, as a $D_V$-bimodule it is isomorphic to $D_V$.
First, let us show that the unit map $N \to D_V^* \tensor \Omega_V \tensor N$ is coacyclic. As in the proof of \ref{prop_filt_A} the cone
of this map is the totalization of the bicomplex
$$
B^{pq} = \begin{cases}
(D_V^* \tensor \Omega^p_V) \tensor N^q,&\text{if $p \ge 0$}\\
N^q,&\text{if $p=-1$}.
\end{cases}
$$
The stages of the canonical filtration of $B$ by vertical slices are total complexes of finite exact sequences of $\Omega_X \tensor D_V$-modules,
and are absolutely acyclic. Therefore, the cone itself is coacyclic.

It remains to show that the counit $\Omega_V \tensor D_V^* \tensor M \to M$ is coacyclic in $\Omega_A$-modules. $D_V^*$ is equipped with
a finite filtration by $\O_V$-subbimodules dual to the standard filtration on $D_V$. It induces a finite filtration $G_\bullet K$ on the cone
of the counit map $K$.

This filtration respects the action of $\Omega$, but not the differential. However, we have $d(G_i) \subset G_{i+1}$, and so it induces a differential
on the associated graded
$$
(\Omega_V \tensor_k \Lambda^\bullet V^*) \tensor_{\O_V} M \to M,
$$
turning it into the standard resolution of $M$ as an $\Omega_V$-module over $\O_V$. Since it is exact we can modify filtration $G$ into
a filtration $\tilde G$ by setting $\tilde G_i = G_i + d(G_i)$, which is closed under differential and has contractible quotients. Therefore,
the cone $K$ is absolutely acyclic.

\qed

\begin{remark}
\label{rem_D_Omega}
\begin{enumerate}[a)]
\item
In the proof we first dualized even variables then odd, however one can proceed the other way around, first dualizing odd variables
and then even without changing the proof.

\item
Since the underlying graded algebra $\# \Omega_A$ is Gorenstein we have equivalences
$$
\D(D_A\Mod) \equ \D^\co(\Omega_A\Mod) \equ \D^\cnt(\Omega_A\Mod).
$$
One can obtain the equivalence of the derived category of $D$-modules and contra-derived category of $\Omega$-modules directly
by using other pairs of adjunctions that behave nicely with respect to limits:
$$
\xymatrix@C=6em{
\D(D_A\Mod) \ar^-{\Omega_X \tensor -}@<0.2em>[r] & \ar^-{\Hom(D_X, -)}@<0.2em>[l] \D^\cnt(\Omega_X \tensor D_V\Mod)
\ar^-{\Hom(\Omega_V, -)}@<0.2em>[r] & \ar^-{D_V \tensor -}@<0.2em>[l] \D^\cnt(\Omega_A\Mod).
}
$$

\item
The equivalence of the theorem is obtained as the composition of two different adjunctions and is not functorial as it depends on the
choice of the trivialization. Notice however, that the functor $\Hom(\Omega_X, -)$ is isomorphic to  $\omega_X^{-1} \tensor \Omega_X \tensor -$
and the functor $\Hom(D_V, -)$ is isomorphic to $\det V \tensor D_V \tensor -$. Twisting as needed we can rewrite the
equivalence in the functorial form
$$
\xymatrix@C=6em{
\D(D_A\Mod) \ar^-{\omega_A^{-1} \tensor \Omega \tensor -}@<0.2em>[r] & \ar^-{D_A \tensor -}@<0.2em>[l] \D^\co(\Omega_A\Mod),
}
$$
where $\omega_A$ denotes the dualizing sheaf on the $\Z$-super-manifold $(X, E)$. 

\end{enumerate}
\end{remark}

\vskip 3em

\subsection{Filtered $D$- and $\Omega$-modules}

The algebra of differential operators $D_A$ is equipped with the standard increasing filtration $F$, and we will write $\tilde D = \Rees (D, F)$.
On the other hand, algebra $\Omega_A$ in addition to cohomological and $s$-gradings has an $h$-grading by the degree of differential forms.
If we equip $\Omega_A$ with $C$-coaction corresponding to the de Rham differential we obtain an algebra in $C$-comodules. As described in
section \ref{duality_filt_A} we can form a Rees algebra $\tilde \Omega_A$ for the decreasing filtration on $\Omega_A$,
and the coderived categories $\D^\co((\Omega_A, C)\MC)$ and $\D^\co(\tilde \Omega_A\Mod)$ are equivalent.

Notice that algebra $\tilde D_A$ no longer has finite cohomological dimension, so we should work with the coderived category of $D$-modules.
With this modification we can upgrade result of theorem \ref{thm_D_Omega} to the filtered case.

\begin{theorem}
\label{thm_D_Omega_filt}
Let $A$ be the ring of functions on a smooth $\Z$-super-manifold, then we have an equivalence
$$
\D^\co(\tilde D_A\Mod) \equ \D^\co(\tilde \Omega_A\Mod).
$$
\end{theorem}
\proof
The proof is the same as in unfiltered case with the $\Hom$-spaces and tensor products taken now over respective $R$-algebras. One thing
we need to upgrade is the statement that the counit map $\tilde D_X \tensor \wtilde{Sp} \tensor \tilde N \to \tilde N$ in $\tilde D$-modules
is not just a quasi-isomorphism,
but a coacyclic map. In fact we will show that it is absolutely acyclic using argument from the end of the proof of theorem \ref{thm_D_Omega}.

Consider filtration $G_\bullet K$ of the cone of the counit map induced by decreasing filtration of $\wtilde{Sp}$ dual to the
tautological filtration of $\tilde \Omega_X$. It is not compatible with the differentials, however $d(G_i) \subset G_{i+1}$,
therefore passing to the associated graded, we obtain the complex
$$
\left(\tilde D_X \tensor_R \Lambda^\bullet \Der(X) \right) \tensor \tilde N \to \tilde N,
$$
which is the standard resolution of $\tilde N$ as a $\tilde D_X$-module over $\O_X \tensor R$. Since it is exact we can modify
filtration $G$ to get a filtration by subcomplexes $\tilde G_i = G_i + d(G_i)$, such that the associated quotients are contractible.
This is a finite filtration and therefore cone $K$ is absolutely acyclic.

\qed

\vskip 1em
\begin{nparagraph}
Next we look at variations of this theorem in two bounded cases. Recall that the Rees algebras $\tilde D$ and $\tilde \Omega$ possess
$h$-grading corresponding to the filtration. Moreover, $\tilde D$ is positively $h$-graded, while $\tilde\Omega$ can be obtained from
an algebra $\Omega$ in $C$-comodules concentrated in the negative $h$-degrees.

For an $h$-graded dg-algebra $\tilde A$ we will write $\Htp(\tilde A\Mod)^{h+}$ for the homotopy category of $\tilde A$-modules
with $h$-degrees bounded from below. Similarly, we can form localizations of this category $\D(\tilde A\Mod)^{h+}$, $\D^\co(\tilde A\Mod)^{h+}$, etc.
We will also need another localization of the homotopy category. We will say that an $\tilde A$-module $M$ is graded-coacyclic if the
associated graded module $\gr M$ is coacyclic as an $\gr_0 \tilde A$-module. It is clear that we have
inclusions
$$
\CoAc(\tilde A\Mod)^{h+} \into \gr\textrm{-}\CoAc(\tilde A\Mod)^{h+} \into \Acyc(\tilde A\Mod)^{h+}.
$$
We will write $\D^{\grco}(\tilde A\Mod)^{h+}$ for
the localization of the homotopy category with respect to morphisms with graded-coacyclic cones.

The composition of functors providing the equivalence of theorem \ref{thm_D_Omega_filt} does not preserve subcategories with bounded
below $h$-grading in general, however we have the following proposition.
\end{nparagraph}

\begin{proposition}
\label{prop_KD_h+}
Let $A$ be the ring of functions on a smooth $\Z$-super-manifold, then the pair of adjoint functors 
\begin{align*}
\Hom(\tilde \Omega_A, -) &\from \Htp(\tilde D_A\Mod)^{h+} \to \Htp(\tilde \Omega_A\Mod)^{h+}\\
\tilde\Omega_A \tensor - &\from \Htp(\tilde \Omega_A\Mod)^{h+} \to \Htp(\tilde D_A\Mod)^{h+},
\end{align*}
induce equivalences
\begin{align*}
\D^\grco(\tilde D_A\Mod)^{h+} &\equ \D^\co(\tilde \Omega_A\Mod)^{h+},\\
\D(\tilde D_A\Mod)^{h+} &\equ \D(\tilde \Omega_A\Mod)^{h+}.
\end{align*}
\end{proposition}
\proof
Since $\tilde D_A$ is concentrated in positive $h$-degrees and $\Omega_A$ in negative $h$-degrees (and therefore $\tilde \Omega_A$ stabilizes
in positive degrees),
the functors in question restrict to the subcategories of modules with bounded below $h$-grading. Moreover, restrictions
on grading imply that $\Hom(\tilde\Omega_A, N) \equ \wtilde{Sp} \tensor N$, so that both functors preserve filtered colimits,
and send coacyclic objects to coacyclic.

\vskip 1em
\textit{First equivalence.}
Let $M \in \tilde\Omega_A\Mod^{h+}$, then $M$ has a bounded from below increasing filtration by $\tilde\Omega_A$-submodules. If every
associated quotient $\gr_i M$ is coacyclic as an $A = \gr_0 \tilde\Omega_A$-module, then $M$ is coacyclic $\tilde\Omega_A$-module.
So the classes of coacyclic and graded-coacyclic module coincide and $\D^\grco(\tilde\Omega_A\Mod)^{h+} \equ \D^\co(\tilde\Omega_A\Mod)^{h+}$.

Let $N$ be a graded-coacyclic object in $\tilde D_A\Mod^{h+}$ and show that $\wtilde{Sp} \tensor N$ is coacyclic. It has an increasing
filtration by $\tilde\Omega_A$-submodules induced by the tautological filtration of $\wtilde{Sp}$. The associated graded for this filtration
are $A$-modules $\wtilde{Sp}^n \tensor N$ and we must show that they are coacyclic. Indeed,
since $N$ is graded-coacyclic $\tilde D_A$-module, it has an increasing filtration by
$A$-submodules with coacyclic quotients and hence it is also coacyclic as an $A$-module.

In order to show that the cone of the unit $\tilde M \to \wtilde{Sp} \tensor \tilde D \tensor \tilde M$ is coacyclic, use the filtration
induced by the tautological filtration on $\tilde M$. It is an increasing filtration bounded from below by assumption on $\tilde M$ and
therefore in order to show that the cone coacyclic it is enough to show coacyclicity of the associated graded module. But since $\tilde \Omega^{\ge 1}$
acts trivially on $\gr M$, this follows from lemma \ref{spencer_coacyclic}.

It remains to show that the counit map $\tilde D \tensor \wtilde{Sp} \tensor N \to N$ is graded-coacyclic, in other words that
the map between associated quotients $\gr_i \left( \tilde D \tensor \wtilde{Sp} \tensor N \right) \to \gr_i N$ is coacyclic in $A$-modules.
For convenience, denote the module on the left by $L_i$. It has another filtration by subcomplexes of $A$-modules
$G_\bullet L_i$ induced by the tautological increasing filtration on $\tilde D \tensor \wtilde{Sp}$. The associated graded for this
filtration is
$$
\gr_\bullet^G L_i = \bigoplus_{p+q = i} \gr_p (\tilde D \tensor \wtilde{Sp}) \tensor \gr_q N.
$$
Now, the zero-term $\gr^G_0 L_i \isom \gr_i N$, while higher terms are finite exact complexes of $A$-modules. This establishes
that the counit map is graded-coacyclic.

\vskip 1em
\textit{Second equivalence.}
To show that the first equivalence descends to the derived categories we just need to show that the adjoint functors preserve
acyclic objects.

Let $N$ be an acyclic object in $\tilde D_A\Mod^{h+}$, and show that $\wtilde{Sp} \tensor N$ is an acyclic $\tilde\Omega_A$-module.
To do that it is enough to consider increasing bounded from below filtration on $\wtilde{Sp} \tensor N$ induced by the filtration
of $\wtilde{Sp}$. The associated quotients are acyclic, since $N$ is an acyclic $A$-module, therefore using spectral sequence
argument we find that $\wtilde{Sp} \tensor N$ is acyclic.

Now let $M$ be an acyclic object in $\tilde\Omega_A\Mod^{h+}$, and show that $\tilde D \tensor M$ is acyclic. Consider filtration
induced by the filtration of $M$, its associated quotients are acyclic, since $\gr_i M$ are acyclic. Since it is a bounded increasing
filtration we conclude as before that $\tilde D \tensor M$ is also acyclic.

\qed

\begin{remark}
In the classical case the algebra $A = \O_X$ has finite cohomological dimension, so that classes of acyclic and coacyclic objects in $A\Mod$
coincide. Therefore, $\D^\co(\tilde \Omega_X\Mod)^{h+} \equ \D(\tilde \Omega_X\Mod)^{h+}$ and
$\D^\grco(\tilde D_X\Mod)^{h+} \equ \D(\tilde D_X\Mod)^{h+}$.
\end{remark}

\vskip 1em
Assume for a moment that the algebra of functions $A$ on a smooth super-manifold is concentrated in negative $s$-grading, in which
case $\Omega_A$ is also negatively $s$-graded, and $s$-grading of $D_A$ is bounded from below ($A$ has only odd negative generators).
We would like to consider full subcategories $\Htp(\tilde D_A\Mod)^{s+}$, $\Htp(\tilde \Omega_A\Mod)^{s+}$
of the corresponding homotopy categories formed by objects with $s$-grading bounded from below, and their localizations.
The pair of adjoint functors $\Hom(\tilde \Omega_A, -)$ and $\tilde D_A \tensor -$ restrict to these subcategories and we have the following statement.

\begin{proposition}
Let $A$ be the ring of functions on a smooth super-manifold concentrated in negative $s$-degrees, then the following localizations are equivalences
$$
\xymatrix@R=0.5em{
\D^\co(\tilde \Omega_X\Mod)^{h+,s+} \ar^-{\equ}[r] & \D(\tilde \Omega_X\Mod)^{h+,s+},\\
\D^\grco(\tilde D_X\Mod)^{h+,s+} \ar^-{\equ}[r] & \D(\tilde D_X\Mod)^{h+,s+}.
}
$$
\end{proposition}
\proof
It is enough to show that the classes of graded-coacyclic and acyclic modules coincide in each case. Let $M$ be an acyclic
object in $\tilde\Omega_A\Mod^{h+,s+}$, we need to show that the associated graded modules $\gr_i M$ for the tautological
filtration are coacyclic $A$-modules. Indeed, consider an increasing filtration by $A$-submodules $G_p \gr_i M = (\gr_i M)_{s\le p}$.
The quotients $G_p / G_{p-1}$ are acyclic modules over $A / A_{s \le 1} \isom \O_X$. Since $\O_X$ is of finite cohomological
dimension we conclude that they are also coacyclic, and hence $M$ is coacyclic as well. The proof on the $D$-module side is
the same.

\qed

\vskip 3em
\subsection{Equivariant $D$- and $\Omega$-modules}

One of the advantages of working on the $\Omega$-module side is that the analogs of notions of weak and strong equivariance
of $D$-modules are defined in a uniform fashion. Let $G$ be an affine $\Z$-super-group, then the algebra of functions $\O_G$ is a commutative Hopf algebra.
The comultiplication $\Delta$ can be extended from $\O_G$ to the de Rham complex $\Omega_G = \O_G \tensor \Sym^\bullet(\g^*[-1])$
by setting $\Delta(x) = x \tensor 1 + 1 \tensor x$ for any $x \in \g^*$. This makes $\Omega_G$ into a commutative Hopf algebra
in $\dgsVect$.

An affine $\Z$-super-manifold $X$ with the left action of $G$ is given by an algebra $\O_X$ in the category of left $\O_G$-comodules
relative to the monoidal structure as described in section \ref{duality_filt_A}. The coaction of $\O_G$ on $\O_X$ induces
a coaction map $a\from \Omega_X \to \Omega_G \tensor \Omega_X$. The latter is determined by its restriction to $\Omega_X^1$,
and corresponds to the derivation $d\circ a\from \O_X \to \O_G \tensor \O_X \to \Omega_G \tensor \Omega_X$. Furthermore,
restricting with respect to the map $\Omega_G \to \O_G$, we obtain the coaction $\Omega_X \to \O_G \tensor \Omega_X$. In other
words $\Omega_X$ is an algebra is $\Omega_G$-comodules and $\O_G$-comodules.

\begin{definition}
\begin{enumerate}[a)]
\item A weakly $G$-equivariant $\Omega_X$-module is an object of $(\Omega_X, \O_G)\MC$.
\item A strongly $G$-equivariant $\Omega_X$-module is an object of $(\Omega_X, \Omega_G)\MC$.
\end{enumerate}
\end{definition}

The relative version of Koszul duality provides equivalences
\begin{align*}
\D^\co((\Omega_X, \O_G)\MC) &\equ \D^\co(\wbar\Cob_{\Omega_X}(\O_G \tensor \Omega_X)),\\
\D^\co((\Omega_X, \Omega_G)\MC) &\equ \D^\co(\wbar\Cob_{\Omega_X}(\Omega_G \tensor \Omega_X)).
\end{align*}

Similarly, we define weakly and strongly $G$-equivariant filtered $\Omega_X$-modules.

\vskip 1em
\begin{nparagraph}
The main reason we want to work with the filtered $D$-modules instead of unfiltered is that they have the ``expected'' behavior.
In particular, let $X$ and a smooth super-manifold and $G$ super-group acting on $X$. Since $D$-modules, as well as $\Omega$-modules
according to equivalence \ref{thm_D_Omega}, don't detect the odd part of either $X$ or $G$, the
action of the odd part of the super-group is ignored. For our purposes this is not acceptable, since this action plays
an important role in Beilinson-Bernstein localization theory (see example in section \ref{ss_example}).

This problem doesn't arise in the classical case because passing to $D$-modules doesn't lose any information about the
action of the group. Indeed, let $X$ and $Y$ be two smooth manifolds, denote $X_\dR$ the de-Rham stack of $X$, i.e., the quotient
of $X$ by the action of the formal completion of the diagonal $\hat\Delta \subset X \times X$. It is easy to see that if
$f\from Y_\dR \to X_\dR$ has a lift $\tilde f\from Y \to X$, then such a lift is unique.

\end{nparagraph}

\begin{nparagraph}[Symplectic quotient.]
Let $(X, E)$ be a smooth $\Z$-super-manifold and $G$ a $\Z$-super-group, acting on it. The cotangent bundle on $(X, E)$ as
a super-manifold is isomorphic to $(T^*X, E \oplus E^*)$, and the induced action of $G$ on $T^*(X, E)$ preserves the natural
symplectic form. Denote by $\mu\from T^*(X, E) \to \g^*$ the moment map for this action.

The (derived) symplectic quotient of $T^*(X, E)$ by the action of $G$ is defined as
$$
\left(T^*(X, E)\ \mathop{\times}\limits^h_{\g^*}\ \{0\} \right) \bigg/ G,
$$
where $\mathop{\times}\limits^h_{\g^*}$ denotes the homotopy fiber product over $\g^*$ with respect to the moment map.

Denote by $K(\g, \O_{\g^*}) = \Sym^\bullet (\g[1]) \tensor \Sym^\bullet (\g)$ the Koszul complex resolution of the skyscraper sheaf at $0 \in \g^*$.
On each open affine piece $U$ of $(X, E)$ the homotopy fiber product can be expressed as $\Spec K(\g, \O_{T^*U})$.

In the case when $X$ is a $G$-torsor over a super-manifold $Y$ the symplectic quotient $T^*X /\!/ G$ is quasi-isomorphic to the cotangent
bundle $T^*Y$.

\end{nparagraph}

\begin{example}
We will illustrate connection between the strongly equivariant $\Omega$-modules and symplectic quotients with the following example.
Consider the odd affine line $G = \A^{0,1}$ with the group operation given by addition. We have $\O_G = k[\epsilon]$, with the degree
$|\epsilon|_s = 1$ and comultiplication $\Delta(\epsilon) = 1 \tensor \epsilon + \epsilon \tensor 1$, and algebra $\Omega_G = k[\epsilon, d\epsilon]$
with coproduct extended from $\O_G$ by setting $\Delta(d\epsilon) = 1 \tensor d\epsilon + d\epsilon \tensor 1$.

Let us describe strongly $\A^{0,1}$-equivariant $\Omega$-modules on a point $\Spec k$, or in other words $\Omega$-modules on the quotient
stack $B\A^{0,1}$. Since the coalgebra $\Omega_G$ is conilpotent we have
$$
\D^\co(\wbar{\Cob(\Omega_G)}\Mod) \equ \D(\wbar{\Cob(\Omega_G)}\Mod).
$$
Now, since derived categories preserve quasi-isomorphism we may replace the cobar complex with the quasi-isomorphic algebra
$B = k[u, du]$, where $\bar u = 1$ and $|u|_s = 1$. The latter, after shear of grading, can be thought of as the de Rham complex
of affine line $\A^1 = \Spec k[x]$, with the action of $\GG_m$ on $\A^1$ by scaling, which corresponds to the original super-grading.
Since $B$ is quasi-isomorphic to $k$, the derived category $\D(B\Mod)$ and hence the category of strongly $\A^{0,1}$-equivariant $\Omega$-modules
on a point is equivalent to $\dgsVect$. This result agrees with the fact that from the point of view of unfiltered $D$-modules the group
$\A^{0,1}$ is indistinguishable from a point.

\vskip 1em
Now we turn to the case of filtered $\Omega$-modules. As before the filtered coalgebra $\tilde\Omega_G$ is conilpotent, so that
the derived and coderived categories over the reduced cobar complex coincide, and we can pass to the derived category of
$\tilde B$-modules, where $\tilde B = \Rees(B, F)$ for the one-step filtration $F_{-1} B = k[u] du \subset F_0 B = B$.

Let us restrict ourselves to the full subcategories formed by object with bounded from below $h$-grading, in which case according
to proposition \ref{prop_KD_h+} we have an equivalence of derived categories of filtered $D$-modules and filtered $\Omega$-modules.
One can see that $\tilde \Omega_G$ is quasi-isomorphic to $k[h, v] / hv$, where $h$ is the variable from the Rees construction and
$v$ corresponds to $d\epsilon$, with $\bar v = 1$, $|v|_s = 1$ and $|v|_h = -1$. The category of modules over this algebra has
a rich structure, in particular one can distinguish between the structure sheaf $\O_{\A^{0,1}}$, corresponding to module
$k[v^{-1}]$ supported at $0$, and the skyscraper $D$-module at $0 \in \A^{0,1}$, corresponding to module $k[h]$.

Now, the derived category $\D(\tilde B\Mod)^{h+}$ can be shown to be equivalent to $\D^\co(\tilde B\Mod)^{h+}$. After change of
grading as before, the latter can be identified with the coderived category of $\GG_m$-equivariant $\tilde\Omega_{\A^1}$-modules,
or equivalently with the equivariant derived category $\D(D_{\A^1}\Mod)^{\GG_m, h+}$. Under this identification
the differential operator $\d_x$ corresponds to the generator of the Koszul algebra $K(\g, \O_{\g^*})$.

\end{example}


\vfill\eject
\section{Representations of Lie super-algebras}

A $\Z$-graded Lie super-algebra is an object $\g \in \sVect$ equipped with a bracket operation $[,]\from \Lambda^2 \g \to \g$,
where the exterior product is understood with respect to the monoidal structure in $\sVect$. The bracket satisfies the Jacobi
relation, i.e., for any $f, g, h \in \g$ we have
$$
[[f, g], h] = [f, [g, h]] - (-1)^{|f|_s |g|_s}[g, [f, h]].
$$
We will only consider finite dimensional Lie super-algebras, in the sense that the total dimension $\sum_i \dim(\g_{s=i})$ is finite.
Denote by $\g_\ev \subset \g$ and $\g_\odd \subset \g$ the subspaces of even and odd $s$-degrees respectively.

The universal enveloping algebra $U\g$ is the quotient of the tensor power $T^\bullet \g$ by the commutation relations
$fg - (-1)^{|f|_s |g|_s} gf = [f, g]$. For any $\g$-module $M$ we can form two complexes, the chain complex
$$
C_\bullet(\g, M) = \bigoplus_n \Sym^n (\g[1]) \tensor M,
$$
and the cochain complex
$$
C^\bullet(\g, M) = \bigoplus_n \Sym^n (\g^*[-1]) \tensor M,
$$
with the differentials induced by the bracket in $\g$ and the action of $\g$ on $M$. We will write $C_\bullet(\g)$ and $C^\bullet(\g)$
for the chain and cochain complexes of the trivial representation of $\g$. It is easy to see that $C^\bullet(\g)$ is an algebra
and $C_\bullet(\g)$ is a coalgebra in $\dgsVect$.

\vskip 1em
\begin{nparagraph}[Lie super-algebra $\gl(V)$.]
Let $V = \bigoplus_i V_{s=i}$ be an object of $\sVect$, we will write $d_i = \dim V_{s=i}$ and $d = \sum_i d_i$.
The algebra of $k$-endomorphisms of $V$ (not necessarily preserving $s$-grading)
will denoted by $\gl(V)$. It is a super-vector space with components
$$
\gl(V)_{s = n}\ =\ \bigoplus_i\ \Hom(V_{s=i}, V_{s=i+n}).
$$
The $\Z$-graded Lie super-algebra structure on $\gl(V)$ is given by the commutator bracket $[f, g] = fg - (-1)^{|f|_s |g|_s} gf$.

Choosing an ordered basis $\{e_1, \ldots e_d\}$ of $V$ (we do not assume compatibility of this order with the ordering of components
$V_{s=i}$), we can identify $\gl(V)$ with the space of $d \times d$ matrices, such that an entry in position $(i, j)$ is of degree
$|e_i|_s - |e_j|_s$. The subalgebra of upper-triangular matrices with respect to this basis is a Borel subalgebra of $\gl(V)$,
and its nilpotent radical consists of upper-triangular matrices with zeroes on the diagonal.

The Lie subalgebra $\sl(V) \subset \gl(V)$ is formed by endomorphisms $f \in \gl(V)$ with zero super-trace
$$
\str f = \sum_{i} (-1)^i \tr (f|_{V_{s=i}}).
$$
\end{nparagraph}

\vskip 3em
\subsection{Derived and coderived categories of $\g$-modules}

Let us discuss relation between derived and coderived categories of $\g$-modules, $C^\bullet(\g)$-modules and $C_\bullet(\g)$-comodules.
This can be considered a generalization of duality described in section \ref{filt_bicomplex}, where we were working with
$\g$-modules for the odd one-dimensional Lie super-algebra $\g = k^{0,1}$ with trivial bracket. In the following propositions
$\g$ is a finite-dimensional $\Z$-graded Lie super-algebra.

\begin{proposition}
The pair of adjoint functors $C_\bullet(\g, -)\from \g\Mod \to C_\bullet(\g)\Comod$ and $U\g \tensor - \from C_\bullet(\g)\Comod \to \g\Mod$
induce equivalence
$$
\D(U\g\Mod) \equ \D^\co(C_\bullet(\g)\Comod).
$$
\end{proposition}
\proof
Clearly coalgebra $C_\bullet(\g)$ is conilpotent, therefore from the general Koszul duality \ref{prop_KD_bar} we find that
the coderived category of $C_\bullet(\g)$-comodules is equivalent to the derived category of $\wbar{\Cob(C_\bullet(\g))}$-modules.
It is straightforward to check that this cobar algebra is quasi-isomorphic to the universal enveloping algebra $U\g$. Since
derived categories respect quasi-isomorphisms, the equivalence descends to $\D(U\g\Mod)$.

\qed

\begin{proposition}
\label{prop_Ug}
We have an equivalence of categories
$$
\D^\co(U\g\Mod) \equ \D^\co(C^\bullet(\g)\Mod).
$$
\end{proposition}
\proof
We proceed as indicated in the remark \ref{rem_D_Omega}(a) by separating even and odd parts of $\g$, and then dualizing first
the odd part and then the even. Since $\g_\odd$ is a $\g_\ev$-module with respect to the adjoint action,
$\g_\odd^*$ is equipped with the coadjoint action $\ad^*$ of $\g_\ev$ and we may form intermediate algebra
$$
B = U\g_\ev \tensor \Sym^\bullet (\g_\odd^*[-1]).
$$
The product is given by the obvious product on the tensor components and the commutation relation $[x, y] = \ad^*_x (y)$,
for any $x \in \g_\ev$ and $y \in \g_\odd^*$, and put the differential to be $0$. Consider also an element $c \in B^2$
defined as the composition
$$
\xymatrix@C=3em{
c\from k \ar^-{\ev^*}[r] & \g_\ev[1] \tensor \g^*_\ev[-1] \ar@{^(->}[r] & \g[1] \tensor \g^*[-1] \ar^-{1\tensor d}[r] & & &
}
$$
$$
\xymatrix@C=3em{
& & \ar^-{1 \tensor d}[r] &\g[1] \tensor \Sym^2 (\g^*[-1]) \ar^-{1 \tensor \mathrm{pr}}[r] & \g[1] \tensor \Sym^2 (\g^*_\odd[-1]).
}
$$
Here $d$ denotes the differential in the cochain complex $C^\bullet(\g)$, and $\mathrm{pr}$ is the projection to the odd component.
This is a central element. Indeed, from the definition it is immediate that $c$ is $\g_\ev$-invariant, and therefore it commutes with
the first factor $U\g_\ev$ of $B$. Using Jacobi identity one can show that it also commutes with the second factor $\Sym^\bullet(\g^*_\odd[-1])$.
In other words $(B, d = 0, c)$ is a curved dg-algebra, i.e. $d^2(x) = [c, x]$ for any $x \in B$.

A curved $B$-module is a $\# B$-module $M$ equipped with the differential $d$ satisfying the usual identities except that
$d^2(m) = cm$ for any $m \in M$. Denote by $B\CMod$ the category of curved modules over the curved dg-algebra $(B, d, c)$, and morphisms
commuting with the differential.

We have the following pair of adjunctions
$$
\xymatrix@C=7em{
\D^\co(U\g\Mod) \ar^-{\Sym^\bullet(\g_\odd^*[-1]) \tensor -}@<0.2em>[r] & \ar^-{\Hom(U\g_\odd, -)}@<0.2em>[l] \D^\co(B\CMod)
\ar^-{\Hom(C^\bullet(\g_\ev), -)}@<0.2em>[r] & \ar^-{U\g_\ev \tensor -}@<0.2em>[l] \D^\co(C^\bullet(\g)\Mod).
}
$$

Using the same argument as in the proof of theorem \ref{thm_D_Omega} one shows that the first adjunction is an equivalence,
furthermore, using the refined argument from the proof of theorem \ref{thm_D_Omega_filt} one shows that the second
adjunction is also an equivalence.

\qed

\begin{remark}
\begin{enumerate}[a)]
\item Since algebra $B$ is of finite cohomological dimension, classes of coacyclic and contra-acyclic curved $B$-modules coincide. However,
the notion of derived category of curved modules is not well-defined, since we can not talk about quasi-isomorphisms.

\item In the classical case the algebra $U\g$ itself has finite cohomological dimension, therefore we have equivalences
$$
\D(U\g\Mod) \equ \D^\co(U\g\Mod) \equ \D^\co(C^\bullet(\g)\Mod).
$$

\item As in the remark \ref{rem_D_Omega}(c) we can twist the equivalence from the proof of the proposition to obtain a functorial
equivalence
$$
\xymatrix@C=7em{
\D^\co(U\g\Mod) \ar^-{(\omega_\g)^{-1} \tensor C^\bullet(\g, -)}@<0.2em>[r] & \ar^-{U\g \tensor -}@<0.2em>[l] \D^\co(C^\bullet(\g)\Mod),
}
$$
where $\omega_\g = \Ber\g[-\dim \g_\ev]$.
\end{enumerate}
\end{remark}

\vskip 1em
\begin{nparagraph}
The universal enveloping algebra $U\g$ is equipped with the standard increasing filtration induced by the tensor powers of $\g$. Similarly
the chain complex $C_\bullet(\g)$ has an increasing filtration and the cochain complex $C^\bullet(\g)$ the dual decreasing filtration.
We will denote by $\wtilde{U\g}$, $\wtilde{C}_\bullet(\g)$ and $\wtilde{C}^\bullet(\g)$ the Rees constructions with respect to these
filtrations.
\end{nparagraph}

\begin{proposition}
\label{prop_Ug_h+}
The pair of adjoint functors
\begin{align*}
\Hom(\wtilde{ C}^\bullet(\g), -) &\from \wtilde {U\g}\Mod \to \wtilde{ C}^\bullet(\g)\Mod,\\
\wtilde {U\g} \tensor - &\from \wtilde {C}^\bullet(\g)\Mod \to \wtilde{ U\g}\Mod,
\end{align*}
induce equivalence
$$
\D(\wtilde {U\g}\Mod)^{h+} \equ \D(\wtilde {C}^\bullet(\g)\Mod)^{h+}.
$$
\end{proposition}
\proof
The proof is similar to proposition \ref{prop_KD_h+}, in fact in this case derived and coderived categories of $\wtilde C^\bullet(\g)$-modules
as well as derived and graded-coderived categories of $\wtilde{U\g}$-modules are equivalent.

\qed

\vskip 3em
\subsection{Localization of $\g$-modules}

Let $G$ be an affine algebraic $\Z$-super-group, and $\g$ its Lie algebra. Consider a smooth $\Z$-super-manifold $X$ with a
transitive action of $G$. We will write $\g_X = \g \tensor \O_X$ for the trivial bundle over $X$ with the fiber $\g$. The action
of $G$ on $X$ induces a map of Lie algebras $a\from \g_X \to TX$, where $TX$ is the tangent bundle of $X$. The kernel of
this map $\s_X \subset \g_X$ is the bundle of stabilizers on $X$.

We will denote by $U\g_X$ the sheaf of $\O_X$-algebras $\O_X \tensor U\g$ with the product determined by the commutation relation
$$
[x, f] = a(x)(f),
$$
for $x \in \g$ and $f \in \O_X$. Thus, we have an $\O_X$-algebra map $a\from U\g_X \to D_X$ respecting the standard filtrations on $U\g$ and
$D_X$.

Furthermore, we will denote by $C^\bullet(\g)_X = \O_X \tensor C^\bullet(\g)$ the commutative dg-algebra with the differential given by
$$
d(f \tensor c) = a^*(df) \tensor c + (-1)^{|f|} f \tensor d_{C^\bullet(\g)}(c),
$$
for $f \in \O_X$, $c \in C^\bullet(\g)$. Here $a^*\from T^*X \to \g^*_X$ is the map dual to $a$. It allows us to connect algebras
$C^\bullet(\g)$ and $\Omega_X$ via the following diagram
$$
\xymatrix{
\Omega_X \ar^-i[r] & C^\bullet(\g)_X & \ar_-p[l] C^\bullet(\g),
}
$$
with both arrows respecting the usual filtrations.

\vskip 1em
\begin{nparagraph}
The forgetful functor $U\from D_X\Mod \to U\g_X\Mod$ clearly sends coacyclic objects to coacyclic, therefore it descends to
the coderived categories. Furthermore, using equivalence $\D^\co(U\g_X\Mod) \equ \Htp(U\g_X\shpInj)$ from proposition \ref{prop_posit}(b)
and similarly for $U\g$-modules we obtain functor $\Gamma\from \D^\co(U\g_X\Mod) \to \D^\co(U\g\Mod)$ by applying the functor of global sections
on $X$ to the injective resolutions.

Similarly, on the $\Omega$-module side we have the global sections functor $\Gamma\from \D^\co(C^\bullet(\g)_X\Mod) \to \D^\co(C^\bullet(\g)\Mod)$.
Consider also functor $F\from \Omega_X\Mod \to C^\bullet(\g)_X\Mod$ defined by
$$
F(M) = (\omega_{\s_X})^{-1} \tensor C^\bullet(\g)_X \tensor_{\Omega_X} M = (\omega_{\s_X})^{-1} \tensor C^\bullet(\s_X, M).
$$
Clearly, $F$ descends to the coderived categories. These functors are compatible with the equivalences of theorem \ref{thm_D_Omega}
and proposition \ref{prop_Ug}.
\end{nparagraph}

\begin{lemma}
\label{lemma_global_sections}
The following diagram commutes.
$$
\xymatrix@C=6em@R=4em{
\D^\co(D_X\Mod) \ar^-{\omega_X^{-1} \tensor \Omega_X \tensor -}_{\equ}[r] \ar_U[d] & \D^\co(\Omega_X\Mod) \ar^{F}[d] \\
\D^\co(U\g_X\Mod) \ar^-{}_{\equ}[r] \ar_{\Gamma}[d] & \D^\co(C^\bullet(\g)_X\Mod) \ar^{\Gamma}[d] \\
\D^\co(U\g\Mod) \ar^-{(\omega_\g)^{-1} \tensor C^\bullet(\g, -)}_{\equ}[r] & \D^\co(C^\bullet(\g)\Mod).
}
$$
\end{lemma}
\proof
The lower square is obviously commutative. For the top square it is easy to see that the composition of $F$ and the top equivalence,
given by $(\omega_{\s_X})^{-1} \tensor \omega_X^{-1} \tensor C^\bullet(\s_X) \tensor \Omega_X \tensor -$, is isomorphic to the other composition
even at the level of $C^\bullet(\g)_X$-modules before passing to the coderived categories.

\qed

\vskip 1em
Functor $F$ has a continuous right adjoint, namely the twisted forgetful functor
$$
\omega_{\s_X} \tensor i_*\from \D^\co(C^\bullet(\g)_X\Mod) \to \D^\co(\Omega_X\Mod).
$$
On the $D$-module side the corresponding functor can be computed by taking injective resolution in $U\g_X$-modules and then applying
$\Hom_{U\g_X}(D_X, -)$. The global sections functor has both left and right adjoints, that will be denoted by $p^*$ and $p^!$ respectively.

\begin{definition}
The (right) localization functor $R^\Omega_X\from \D^\co(C^\bullet(\g)\Mod) \to \D^\co(\Omega_X\Mod)$ is given by the composition
$R^\Omega_X (M) = \omega_{\s_X} \tensor i_* p^! (M)$. On the $D$-module side the corresponding functor is
$R^D_X (V) = \Hom_{U\g_X}(D_X, p^! V)$.
\end{definition}

\begin{remark}
\label{rem_R_loc}
\begin{enumerate}[a)]
\item The lemma holds for filtered modules as well and we define localization functors $\tilde R^\Omega_X$ and $\tilde R^D_X$ in the filtered context
same as above.

\item If $X$ is a compact super-manifold, i.e. the underlying even manifold is compact, the right adjoint $p^!$ is continuous, so that
the localization functor $R^\Omega_X$ is also continuous.

\item Functor $F$ in general doesn't have a left adjoint, this is the case when the stabilizer bundle $\s_X$ has a non-trivial odd component.
Similarly, functor $U$ restricted to the coderived categories doesn't preserve limits in general, since the localization functor
from a homotopy category to its coderived category doesn't preserve them.

\item In the classical case, however, functor $F$ and therefore $U$ as well have both left and right adjoints. Notice that $F(M)$ can
be written as $\Hom_{\Omega_X}(C^\bullet(\g)_X, M)$, hence its left adjoint is the forgetful functor $i_*$. We define the left
localization functor $L^\Omega_X\from \D^\co(C^\bullet(\g)\Mod) \to \D^\co(\Omega_X\Mod)$ by
$$
L_X^\Omega(M) = i_*\left(C^\bullet(\g)_X \tensor_{C^\bullet(\g)} M \right).
$$
On the $D$-module side the coderived categories are equivalent to the derived categories and the left localization $L^D_X$ is given by
$$
L^D_X(V) = D_X \tensor^L_{U\g} V.
$$

\end{enumerate}
\end{remark}

\vskip 1em
\begin{nparagraph}[Localization in the bounded filtered case.]
Let us discuss the localization in the case of filtered modules with bounded from below filtration. Clearly, both functors $U$ and
$\Gamma$ restrict to the corresponding subcategories. On the $\Omega$-module side consider functor
$F'\from \D^\co(\tilde\Omega_X\Mod)^{h+} \to \D^\co(\tilde C^\bullet(\g)_X\Mod)^{h+}$:
$$
F'(M) = \Hom_{\tilde\Omega_X} (\tilde C^\bullet(\g)_X, M) = \tilde C_\bullet(\s_X, M).
$$

As before it is straightforward to check that these functors are compatible with equivalences of propositions \ref{prop_KD_h+} and
\ref{prop_Ug_h+}, in the sense that the following diagram commutes.
$$
\xymatrix@C=6em@R=4em{
\D(\tilde D_X\Mod)^{h+} \ar^-{\Hom(\tilde\Omega_X,  -)}_{\equ}[r] \ar_I[d] & \D(\tilde\Omega_X\Mod)^{h+} \ar^{I}[d] \\
\D^\grco(\tilde D_X\Mod)^{h+} \ar^-{\Hom(\tilde\Omega_X,  -)}_{\equ}[r] \ar_U[d] & \D^\co(\tilde\Omega_X\Mod)^{h+} \ar^{F'}[d] \\
\D^\grco(\tilde U\g_X\Mod)^{h+} \ar^-{}_{\equ}[r] \ar_{\Gamma}[d] & \D^\co(\tilde C^\bullet(\g)_X\Mod)^{h+} \ar^{\Gamma}[d] \\
\D(\tilde U\g\Mod)^{h+} \ar^-{\Hom(\tilde C^\bullet(\g), -)}_{\equ}[r] & \D(\tilde C^\bullet(\g)\Mod)^{h+}.
}
$$
Here $I$ denotes the inclusion functor, right adjoint to the localization with respect to quasi-isomorphisms, in other words it is the
injective resolution functor.

In this case $F'$ has left adjoint given by the forgetful functor $i_*\from \tilde C^\bullet(\g)_X\Mod \to \tilde\Omega_X\Mod$.
We define the left localization $\tilde L^\Omega_X$ as the composition
$$
\tilde L^\Omega_X(M) = i_* \left( \tilde C^\bullet(\g)_X \tensor_{\tilde C^\bullet(\g)} M \right).
$$
On the $D$-module side, passing further to the derived category of $\tilde D_X$-modules we have as in the classical case
$\tilde L^D_X(V) = \tilde D_X \tensor^{L}_{\tilde U\g} V \in \D(\tilde D_X\Mod)$.

\end{nparagraph}

\vskip 1em
\begin{nparagraph}[Monodromic localization.]
Consider two smooth $\Z$-super-manifolds $X$ and $Y$ with transitive actions of $G$, and a $G$-invariant map $f\from Y \to X$.
Assume that $f$ makes $Y$ into a right $K$-torsor over $X$ for some affine algebraic super-group $K$, so that actions of
$K$ and $G$ commute. For a sheaf $E$ on $Y$ we will denote $E^K$ the sheaf on $X$ of $K$-invariants of the direct image $f_* E$.
We will call $\Omega_Y^K$-modules and $D_Y^K$-modules the monodromic $\Omega$-, and $D$-modules on $X$ respectively.

Using the diagram of algebras $\xymatrix@1{\Omega_X \ar[r] & \Omega_Y^K \ar^-i[r] & C^\bullet(\g)_Y^K \isom C^\bullet(\g)_X & \ar_-p[l] C^\bullet(\g)}$
we construct the functor $F$ as before
by putting
$$
F(M) = (\omega_{\s_Y^K})^{-1} \tensor C^\bullet(\g)_X \tensor_{\Omega_Y^K} M \isom (\omega_{\s_Y^K})^{-1} \tensor C^\bullet(\s_Y^K, M).
$$
We define the right localization functor as the right adjoint to the composition $\Gamma F$, i.e.
$R^{\Omega,K}_Y(M) = \omega_{\s_Y^K} \tensor i_* p^!(M)$. On the $D$-module side the corresponding functor is given by
$R^{D,K}_Y(V) = \Hom_{U\g_X}(D_Y^K, p^! V)$.

\end{nparagraph}

\vskip 3em
\subsection{Monadicity of localization}
Consider two dg-categories $\C$ and $\cD$ and a pair of adjoint functors $(L, R)$ between them, specifically $L\from \C \to \cD$ and
$R\from \cD \to \C$. The composition $A = RL\from \C \to \C$ is a monad, and we will write $A\Mod_\C$ for the category of $A$-modules
in $\C$, i.e. collections $(M \in \C, a_1\from A(M) \to M, a_2\from A^2(M) \to M, \ldots)$,
satisfying higher associativity conditions and compatibility with the monad structure of $A$.

The right adjoint $R$ induces the comparison functor $K$.
$$
\xymatrix@R=4em{
\cD \ar@<0.2em>^R[dr] \ar^K[rr] && \ar[dl] A\Mod_\C \\
& \ar@<0.2em>^L[ul] \C &
}
$$
Functor $R$ is called {\em monadic} if the comparison functor is an equivalence.

If $\cD$ has all colimits then $K$ has a left adjoint $J\from A\Mod_\C \to \cD$. Indeed, consider an $A$-module $(M, a)$, and form
the following simplicial object in $\cD$, where the face maps are given by the counit of adjunction $\eta\from LR \to \Id$ and
the action map $a\from RL(M) \to M$:
$$
\xymatrix@C=4em{
\ldots LRLRL (M) \ar@<0.4em>^-{\eta}[r] \ar[r] \ar@<-0.4em>_-{L(a)}[r]  & LRL(M) \ar@<0.2em>^-{\eta}[r] \ar@<-0.2em>_-{L(a)}[r] & L(M).
}
$$
It is immediate to see that the realization (or the homotopy colimit) of this simplicial object gives the left adjoint $J(M)$.
Furthermore, if $R$ preserves colimits then the natural map $M \to KJ(M)$ is an isomorphism, or equivalently, $J$ is a fully
faithful embedding $A\Mod_\C \into \cD$. If in addition $R$ is conservative then the pair $(J, K)$ is an equivalence of
categories, and $R$ is monadic.

\begin{nparagraph}[Ind-monadicity.]
Let $X_\bullet$ be a simplicial object in $\cD$. We say that $X_\bullet$ is of {\em compact type} (relative to functor
$L\from \C \to \cD$) if all $X_i \isom L(Q_i)$ for
some compact objects $Q_i \in \C$, and all but finite number of homotopy groups $\pi_i(X_\bullet)$ vanish. We also say that
an object $Z \in \cD$ is {\em compactly generated} (relative to $L\from \C \to \cD$) if it is isomorphic to the homotopy colimit of
some simplicial object $X_\bullet$ of compact type. We will denote $\cD^\cg$ the full subcategory of compactly generated objects in $\cD$.

\begin{definition}
We say that functor $R$ is {\em ind-monadic} if the induced comparison functor $\hat K\from \cD \to \Ind (A\Mod_\C^\cg)$ is
an equivalence.
\end{definition}
\end{nparagraph}

Let us give a sufficient condition for a functor to be ind-monadic that is useful in the localization setting.

\begin{proposition}
\label{prop_ind_monadic}
Let $(L, R)$ be an adjunction between two compactly generated complete dg-categories $\C$ and $\cD$, such that
\begin{enumerate}[a)]
\item functor $R$ preserves filtered colimits and is conservative on compacts,
\item the full subcategory $\cD_c$ of compact objects in $\cD$ is contained in $\cD^\cg$,
\item if $X_\bullet$ is a simplicial object with $X_i \in \cD_c$, then all $\pi_i(X_\bullet) \in \cD_c$.
\end{enumerate}
Then $R$ is ind-monadic.
\end{proposition}
\proof
First of all notice that since $R$ preserves filtered colimits the left adjoint $L$ sends compact objects in $\C$ to compact
objects in $\cD$. Now, let $Z = \hocolim X_\bullet$ an object in $\cD^\cg$ and $X_\bullet$ a simplicial object of compact type.
We have all $X_i \in \cD_c$, so from (c) we have all $\pi_i(X)$ are compact and since only finite number of them are non-zero
the realization $Z$ is also compact. We have shown that $\cD^\cg$ is contained in $\cD_c$, and combining this with (b) we conclude
that $\cD_c = \cD^\cg$.

Next let us show the following fully faithful inclusions
$$
\xymatrix@C=3em{
\cD_c \ar@{^{(}->}^-K[r] & A\Mod_{\C} \ar@{^{(}->}^-J[r] & \cD.
}
$$
Indeed, a continuous exact functor $R$ between dg-categories commutes with arbitrary colimits, and therefore, as was mentioned in the beginning of this
section, we have the second inclusion. Furthermore, since $K$ is conservative on compacts it identifies $\cD_c$ with a full subcategory of $A$-modules.

It remains to show that $K$ identifies $\cD_c$ with $A\Mod_\C^\cg$. Since $K$ preserves arbitrary colimits, it respects
Postnikov tower construction, hence if $X_\bullet$ has only finite number of non-zero homotopy groups, so does $K(X_\bullet)$.
By definition $KL(Q)$ is the free $A$-module generated by $Q$, therefore $K$ sends simplicial objects of compact type in $\cD$ to
objects of compact type in $A\Mod_\C$ and compactly generated objects in $\cD$ to $A\Mod_\C^\cg$, which establishes one inclusion.

Let $Z = \hocolim X_\bullet$ be a compactly generated object in
$A\Mod_\C$. Since $Z \isom KJ(Z)$, it is enough to show that $J(Z)$ is compactly generated object in $\cD$. $J$ is a left adjoint to $K$,
therefore it preserves arbitrary colimits, so it remains to show that $J(X_\bullet)$ is a simplicial object of compact type in $\cD$.
This follows from the fact that $JKL(Q) \isom L(Q)$ for a compact $Q \in \C_c$.

\qed

\vskip 1em
\begin{nparagraph}[Localization on the flag variety.]
Now, let $G = GL(V)$ for some super-vector space $V$ and $\g = \gl(V)$ its Lie super-algebra. Consider a Borel subgroup $B \subset G$
and its unipotent radical $N \subset B$ for some choice of ordered basis of $V$, as well as the corresponding Lie algebras $\nn \subset \bb \subset \g$.
The flag super-manifold $X$ is the
quotient $G / B$, and the affine flag $\tilde X$ is the quotient $G / N$. The projection map $\tilde X \to X$ is a right $T$-torsor,
where $T$ is the torus of diagonal matrices in $GL(V)$ with respect to the chosen ordered basis. The stabilizer bundle $\s_X$ is
isomorphic to $\bb_X$, with the fibers conjugate to the Borel subalgebra $\bb$, and the bundle $\s^T_{\tilde X}$ is isomorphic
to $\nn_X$ with fibers conjugate to subalgebra $\nn$.

We will be interested in properties
of the right monodromic localization functor $R^{D,T}_{\tilde X}$.
\end{nparagraph}

\begin{lemma}
\label{lemma_conserv_compacts}
The right localization functor $R^{D,T}_{\tilde X}\from \D^\co(U\g\Mod) \to \D^\co(D^T_{\tilde X}\Mod)$ is conservative on compacts.
\end{lemma}
\proof
Recall that coderived category of $U\g$-modules is the ind-completion of the absolute derived category of finitely generated
$U\g$-modules (\ref{prop_compact_gen}), and since $U\g$ sits in cohomological degree $0$, the latter is equivalent to the
bounded derived category of finitely generated modules. Consider the following diagram
$$
\xymatrix@=3em{
\D^\b(U\g\Mod_\fg) \ar_{R^{D,T}_{\tilde{X}}}[d] \ar[r] & \D^\b(U\g_\ev\Mod_\fg) \ar^{R^{D,T}_{\tilde{X}_\ev}}[d] \\
\D(D_{\tilde{X}}^T\Mod) \ar[r] & \D(D_{\tilde{X}_\ev}^T\Mod).
}
$$
The top arrow, given by the restriction of action from $\g$ to $\g_\ev$ is conservative, since quasi-isomorphisms are
detected on the underlying vector spaces. The right arrow is conservative by the classical monodromic localization
result. Therefore, the left arrow is also conservative.

\qed

Denote by $W$ the {\em Weyl monad}, i.e. the endofunctor on the category of monodromic $D_X$-modules obtained as the
composition of global sections and localization functors:
$$
W = R^{D, T}_{\tilde X} \Gamma \from \D^\co(D^T_{\tilde X}\Mod) \to \D^\co(D^T_{\tilde X}\Mod).
$$




\vskip 1em
\begin{theorem}
\label{thm_ind_monadic}
The right localization functor $R^{D, T}_{\tilde X}$ is ind-monadic, in other words we have an equivalence
$$
\xymatrix@C=3em{
\hat K\from\D^\co(U\g\Mod) \ar^-{\equ}[r] & \Ind(W\Mod^\cg_{\D^\co(D^T_{\tilde X}\Mod)}).
}
$$
\end{theorem}
\proof
Let us check conditions of the proposition \ref{prop_ind_monadic}. As was mentioned earlier (\ref{rem_R_loc}(b)) the right localization
functor is continuous, combined with lemma \ref{lemma_conserv_compacts} it establishes condition (a).

To check (c) recall that the coderived category of $U\g$-modules is compactly generated by finitely generated $U\g$-modules. Since
$U\g$ is Noetherian, a complex of finitely generated modules has finitely generated cohomology groups.

It is easier to verify (b) on the $\Omega$-module side. Again, since $C^\bullet(\g)$ is Noetherian, a finitely generated $C^\bullet(g)$-module
has a simplicial resolution of compact type relative to vector spaces. However, from the description of the global sections functor
(lemma \ref{lemma_global_sections}) we immediately see that it is also of compact type relative to monodromic $\Omega$-modules.

\qed

\vskip 3em
\subsection{Example: $\sl(1, 1)$}
\label{ss_example}
Let $V = k \oplus k\<-1\>$ and consider $\g = \sl(1, 1) \subset \gl(V)$, i.e. the Lie subalgebra consisting of elements with
super-trace zero. Fixing an ordered basis of $V$ we identify $\gl(V)$ with a Lie algebra of $2\times 2$ matrices and $\g$ is the subalgebra
spanned by
$$
e = \begin{pmatrix}0& 1\\ 0& 0\end{pmatrix},\quad\quad
f = \begin{pmatrix}0& 0\\ 1& 0\end{pmatrix},\quad\quad
h = \begin{pmatrix}1& 0\\ 0& 1\end{pmatrix}.
$$
The super-degrees are $|e|_s = 1$, $|f|_s = -1$, $|h|_s = 0$, and the bracket is given by
$$
[e, f] = h, \quad\quad [h, e] = [h, f] = 0.
$$
The universal enveloping algebra $U\g$ is isomorphic to a quotient of the free associative algebra on two generators $e$ and $f$ by the following relations
$$
U\g = k\<e, f\> / e^2 = f^2 = 0.
$$

The corresponding group $G = \SL(1, 1)$ is an affine Lie super-group with the ring of functions $\O_G = k[x_{11}, x_{12}, x_{21}, x_{22}] / \Ber = 1$,
where
$$
\Ber = {x_{11} \over x_{22} } - { x_{12} x_{21} \over x_{22}^2},
$$
and coproduct is given by $\Delta(x_{ij}) = \sum_k x_{ik} \tensor x_{kj}$.

\vskip 1em
\begin{nparagraph}[Affine flag variety.]
The affine flag variety $\tilde X = G/N$ is isomorphic to $\GG_m \times \A^{0,1} = \Spec k[t^{\pm}, \epsilon]$, where coordinates $t$ and $\epsilon$
lift to the right $N$-invariant functions $x_{11}$ and $x_{21}$ respectively. The left action of $G$ on $\tilde X$ induces a map
$\g \to D_{\tilde X}$ given by
$$
e \mapsto \epsilon \d_t,\quad\quad f \mapsto t\d_\epsilon,\quad\quad h \mapsto t\d_t + \epsilon\d_\epsilon.
$$

Both functions $t$ and $\epsilon$ are of weight $1$ for the action of the diagonal torus $T$, therefore the space of $T$-invariants
$D^T_{\tilde X}$ is generated as a subalgebra of $D_{\tilde X}$ by $\epsilon t^{-1}$, $t\d_t$ and $t\d_\epsilon$.
Alternatively, using the map $U\g \to D^T_{\tilde X}$ described above, we can rewrite it as
\begin{equation}
\label{eqn_example}
D^T_{\tilde X}\ \isom\ k[h] \tensor k\<\alpha, f\> / [\alpha, f] = 1,
\end{equation}
or equivalently,
$$
D^T_{\tilde X}\ \isom\ U\g\<\alpha\> \big/ e\alpha = \alpha e = 0, [f, \alpha] = 1,
$$
where $\alpha$ corresponds to the element $\epsilon t^{-1}$.

Passing to the associated graded we find that the ring of functions on the bundle $\bb^*_X = (T^* \tilde X) / T$ is the free commutative algebra
$$
\O_{\bb^*_X} = k[\alpha, v_\alpha, h].
$$
The moment map $\mu^*\from\O_{\g^*} \to \O_{\bb^*_X}$ is given by
$$
e \mapsto \alpha h,\quad\quad f \mapsto v_\alpha,\quad\quad h \mapsto h,
$$
in other words
$$
\O_{\bb^*_X}\ \isom\ \O_{\g^*}[\alpha] / e = \alpha h.
$$
Geometrically, the moment map $\mu\from \bb^*_X \to \g^*$ is a contraction of the odd affine line fibers $\Spec k[\alpha]$
over the closed subscheme $Z_{h = 0} \subset \g^*$, and isomorphism on the complement.

\end{nparagraph}

\vskip 1em
\begin{nparagraph}[Weyl monad.]
The second factor in \ref{eqn_example} is an algebra of differential operators on an odd affine line, and as was discussed in the
beginning of section \ref{ss_D_Omega} the category of $D$-modules on an odd affine line is equivalent to $\sVect$. Therefore, the category
of $D^T_{\tilde X}$-modules is equivalent to $k[h]$-modules in $\sVect$. The Weyl monad is represented by an algebra object in
$D^T_{\tilde X}$-bimodules, and since in this case the flag manifold is affine, this algebra can be written as
$$
W = \End_{U\g}(D^T_{\tilde X}).
$$
The enveloping algebra of $D^T_{\tilde X}$ is equivalent to the tensor product of $k[h_1, h_2]$ and the algebra of differential
operators on $\A^{0, 2}$, therefore we can write
$$
W = \wbar W \tensor \O_{\A^{0, 2}},
$$
where $\wbar W$ is some $k[h_1, h_2]$-algebra and $\O_{\A^{0, 2}}$ is the unique indecomposable $D_{\A^{0,2}}$-module.
Moreover, since $h$ is central in $D^T_{\tilde X}$ it is also central in $W$, and $\wbar W$ is in fact a $k[h]$-algebra, where $h = h_1 = h_2$.

The standard resolution
$$
\Sym^\bullet(\nn[1]) \tensor U\g_X \to D^T_{\tilde X}
$$
is periodic and can be written explicitly as
$$
\xymatrix@C=6em{
\ldots U\g^2 \ar^{\begin{pmatrix}e& 0\\0& e\end{pmatrix}}[r] &
U\g^2 \ar^{\begin{pmatrix}0& e\\e& -ef\end{pmatrix}}[r] &
U\g^2 \ar^{\begin{pmatrix}1\\\alpha\end{pmatrix}}[r] & D^T_{\tilde X}.
}
$$
A straightforward calculation shows that
$$
\wbar W \isom k[h, u] / hu,\quad\quad \bar u = 1,\quad |u|_s = -1,
$$
where $u$ is the variable Koszul dual to $e$.

\end{nparagraph}

\vskip 1em
\begin{nparagraph}[Localization in weight $h \neq 0$.]
The category of representations of $\sl(1, 1)$ in non-zero weight $h = \lambda$ is semi-simple with a unique indecomposable $V_\lambda$ up to shift of
super-grading. On the other hand, the quotient $D^T_{\tilde X} / (h - \lambda)$ is isomorphic after rescaling to $D_X$, and so the category of
monodromic $D$-modules on $X$ in weight $h = \lambda$ is equivalent to $D_X$-modules. Clearly under this identification
the unique indecomposable representation $V_\lambda$ localizes as the unique indecomposable $D_X$-module $\O_X$.

\end{nparagraph}

\vskip 1em
\begin{nparagraph}[Localization in weight $h = 0$.]
Let us clarify what happens in weight $0$. The category of representation is no longer semi-simple, but it is generated by the unique
irreducible module, namely the trivial representation. The subcategory of $\D^\co(U\g\Mod)$ with nilpotent action of $h$ is identified
under the Koszul equivalence \ref{prop_Ug} with the $\D(C^\bullet(\g)\Mod)$ considered as a subcategory of the coderived category.

Now compare localizations for the two opposite Borel subalgebras $B^+$ and $B^-$, namely the upper-triangular matrices and lower-triangular matrices.
The respective Weyl monads are represented by algebras $W^+ = k[h, u] / hu$ and $W^- = k[h, v] / hv$, where $v$ is Koszul dual to $f$.
The singularity categories of $W^+$ and $W^-$ modules are generated by $Y^+ = k[u]$ and $Y^- = k[v]$ respectively. They are the localizations
of the trivial representation on $X^+$ and $X^-$. The subcategories generated
by these objects are equivalent to the derived categories of their endomorphism algebras. We find that
$$
\End_{W^+}(Y^+) = k[u, s] / us\ \isom\ C^\bullet(\g)\ \isom\ k[v, t] / vt = \End_{W^-}(Y^-),
$$
where we identify $s = v$ and $u = t$.

\end{nparagraph}

\vfill\eject

\end{document}